\documentclass[a4paper]{amsproc}

\usepackage{amssymb}
\usepackage{pstricks}
\usepackage{pstcol}
\usepackage{graphicx}
\usepackage{amsmath}
\usepackage{amsfonts}
\usepackage{latexsym}

\theoremstyle{plain}

\theoremstyle{definition}

\numberwithin{equation}{section}

\renewcommand{\geq}{\geqslant}

\title[An elementary proof]{The integrals in Gradshteyn and Ryzhik. \\
Part 14: An elementary evaluation of entry $\mathbf{3.411.5}$}
\subjclass[2000]{Primary 33}


\author[T. Amdeberhan]{Tewodros Amdeberhan}
\address{Department of Mathematics,
Tulane University, New Orleans, LA 70118}
\email{tamdeberhan@math.tulane.edu}

\author[V. Moll]{Victor H. Moll}
\address{Department of Mathematics,
Tulane University, New Orleans, LA 70118}
\email{vhm@math.tulane.edu}

\address{\hfill{\it Received ??, revised ??}
\newline Departamento de Matem\'atica
\newline
Universidad T\'ecnica Federico Santa Mar\'{\i}a
\newline  Casilla 110-V,
\newline Valpara\'{\i}so, Chile}

\thanks{The author wishes to acknowledge the partial support of  
NSF-DMS 0713836.}

\begin{document}

{\begin{flushleft}\baselineskip9pt\scriptsize {\bf SCIENTIA}\newline
Series A: {\it Mathematical Sciences}, Vol. ?? (2010), ??
\newline Universidad T\'ecnica Federico Santa Mar{\'\i}a
\newline Valpara{\'\i}so, Chile
\newline ISSN 0716-8446
\newline {\copyright\space Universidad T\'ecnica Federico Santa
Mar{\'\i}a\space 2010}
\end{flushleft}}

\vspace{10mm} \setcounter{page}{1} \thispagestyle{empty}

\begin{abstract}
An elementary proof of an entry in the table of integrals by 
Gradshteyn and Rhyzik s presented.
\end{abstract}

\maketitle

\newcommand{\nn}{\nonumber}
\newcommand{\ba}{\begin{eqnarray}}
\newcommand{\ea}{\end{eqnarray}}
\newcommand{\ift}{\int_{0}^{\infty}}
\newcommand{\ione}{\int_{0}^{1}}
\newcommand{\ifft}{\int_{- \infty}^{\infty}}
\newcommand{\no}{\noindent}
\newcommand{\realpart}{\mathop{\rm Re}\nolimits}
\newcommand{\imagpart}{\mathop{\rm Im}\nolimits}

\newtheorem{Definition}{\bf Definition}[section]
\newtheorem{Thm}[Definition]{\bf Theorem} 
\newtheorem{Example}[Definition]{\bf Example} 
\newtheorem{Lem}[Definition]{\bf Lemma} 
\newtheorem{Note}[Definition]{\bf Note} 
\newtheorem{Cor}[Definition]{\bf Corollary} 
\newtheorem{Prop}[Definition]{\bf Proposition} 
\newtheorem{Problem}[Definition]{\bf Problem} 
\numberwithin{equation}{section}

\maketitle

\section{Introduction} \label{sec-intro}
\setcounter{equation}{0}

The compilation by I. S. Gradshteyn and  I. M. Ryzhik \cite{gr}
contains about $600$ pages of definite integrals. Some of them are quite 
elementary; for instance, $4.291.1$ 
\begin{equation}
\int_{0}^{1} \frac{\ln(1+x)}{x} \, dx  = \frac{\pi^{2}}{12}
\label{simple-1}
\end{equation}
\noindent
is obtained by expanding the integrand as a power series and using 
the value 
\begin{equation}
\sum_{k=1}^{\infty} \frac{(-1)^{k-1}}{k^{2}} = \frac{\pi^{2}}{12}.
\label{alter}
\end{equation}
\noindent
The latter is reminiscent of the series
\begin{equation}
\sum_{k=1}^{\infty} \frac{1}{k^{2}} = \frac{\pi^{2}}{6}.
\label{zeta2}
\end{equation}
\noindent
The reader will find in \cite{irrbook} many proofs of the 
classical evaluation (\ref{zeta2}). 

Most entries in \cite{gr} appear  quite formidable and their
evaluation requires a variety of methods and ingenuity. Entry 
$4.229.7$
\begin{equation}
\int_{\pi/4}^{\pi/2} \ln \ln \tan x \, dx = 
\frac{\pi}{2} \ln \left( \frac{\Gamma(\tfrac{3}{4}) }{\Gamma(\tfrac{1}{4}) }
\sqrt{2 \pi} \right)
\label{vardi-1}
\end{equation}
\noindent
illustrates this point. Vardi \cite{vardi1} describes 
a good amount of mathematics involved in evaluating (\ref{vardi-1}). 
The integral
is first interpreted in terms of the derivative of the $L$-function
\begin{equation}
L(s) = 1 - \frac{1}{3^{s}} + \frac{1}{5^{s}} - \frac{1}{7^{s}} - \cdots
\end{equation}
\noindent
as
\begin{equation}
\int_{\pi/4}^{\pi/2} \ln \ln \tan x \, dx = 
-\frac{\pi \gamma}{4} + L'(1). 
\end{equation}
\noindent
Here $\gamma$ is {\em Euler's constant}. Then 
$L'(1)$ is computed in terms of the {\em gamma function}. This is an unexpected
procedure.

Any treatise such as \cite{gr}, containing 
large amount of information is bound to have some errors. Some of 
them are easy to fix. For instance, formula $3.511.8$ in the sixth 
edition \cite{gr4} reads
\begin{equation}
\int_{0}^{\infty} \frac{dx}{\cosh^{2}x} = \sqrt{\pi} 
\sum_{k=0}^{\infty} \frac{(-1)^{k}}{\sqrt{2k+1}}. 
\label{form-wrong-1}
\end{equation}
\noindent
The source given for this integral is formula $\mathbf{BI}(98)(25)$ from
the table by Bierens de Haan \cite{bierens1}, where it 
appears as
\begin{equation}
\int_{0}^{\infty} \frac{1}{e^{t}+e^{-t}} \frac{dt}{\sqrt{t}} = 
\sqrt{\pi} \sum_{k=0}^{\infty} \frac{(-1)^{k}}{\sqrt{2k+1}}. 
\end{equation}
\noindent
The change of variable $t = \sqrt{x}$ yields a correct version of
(\ref{form-wrong-1}):
\begin{equation}
\int_{0}^{\infty} \frac{dx}{\cosh(x^{2})} = \sqrt{\pi} 
\sum_{k=0}^{\infty} \frac{(-1)^{k}}{\sqrt{2k+1}}. 
\label{new-version-2}
\end{equation}
\noindent
It is now clear what happened: a typo produced (\ref{form-wrong-1}). In the 
latest edition of the integral table \cite{gr}, the 
editors have decided to replace this 
entry with 
\begin{equation}
\int_{0}^{\infty} \frac{dx}{\cosh^{2}x} = 1. 
\label{new-version-1}
\end{equation}
\noindent
The right-hand side of (\ref{form-wrong-1}) has been corrected.

The advent of computer algebra packages has not made these tables obsolete. 
The latest version of \texttt{Mathematica} evaluates (\ref{new-version-1}) directly, 
but it is unable to produce (\ref{new-version-2}). 

Most of the errors in \cite{gr} are of the type: some parameter has been
mistyped, an exponent has been misplaced, parameters are mistaken to be 
identical (a common mishap is $\mu$ and 
$u$ appearing in the same formula).  Despite of this fact, it is
a remarkable achievement for such an endevour. The accuracy of 
\cite{gr} comes from the 
effort of many generations, beginning with \cite{lindman1} and also
including \cite{klerer1, sheldon1}. 

A different type of error was found by one of the authors. It turns 
out that entry $3.248.5$ of \cite{gr} is incorrect. To explain  the reason
for looking at any specific entry  requires some background. 
The stated entry $3.248.5$ involves the rational function
\begin{equation}
\varphi(x) = 1 + \frac{4x^{2}}{3 (1+x^{2})^{2}}
\label{varphi}
\end{equation}
\noindent
and the result says
\begin{equation}
\int_{0}^{\infty} \frac{dx}{(1+x^{2})^{3/2} \left[ \varphi(x) + 
\sqrt{\varphi(x)} \right]^{1/2}} = \frac{\pi}{2 \
\sqrt{6}}.
\label{wrong-3}
\end{equation}

The encounter begins with the evaluation of 
\begin{equation}
N_{0,4}(a;m) = \int_{0}^{\infty} \frac{dx}{(x^{4}+2ax^{2}+1)^{m+1}} 
\label{quartic-eva}
\end{equation}
\noindent
in the form
\begin{equation}
N_{0,4}(a;m) = \frac{\pi}{2} \frac{P_{m}(a)}{[2(a+1)]^{m+1/2}}, 
\end{equation}
\noindent
where $P_{m}(a)$ is a polynomial of degree $m$. The reader will find in 
\cite{amram, manna-moll-survey} details about (\ref{quartic-eva}) 
and properties 
of the coefficients of $P_{m}$.  It is rather interesting 
that $N_{0,4}(a,m)$ appears in the 
expansion of the double square root function
\begin{equation}
\sqrt{a+ \sqrt{1+ c}} = \sqrt{a+1} + \frac{1}{\pi \sqrt{2}} 
\sum_{k=1}^{\infty} \frac{(-1)^{k-1}}{k} N_{0,4}(a;k-1)c^{k}. 
\end{equation}
\noindent
Browsing \cite{gr4} on a leisure  day, and with double square roots in our 
mind, formula 
$3.248.5$ caught our attention. After many failed attempts at the 
proof, a simple numerical integration showed that (\ref{wrong-3}) is incorrect. 
In spite of our inability to evaluate this integral, we have produced many
equivalent versions. The reader is invited to verify that, 
if $\sigma(x,p) := \sqrt{x^{4} + 2px^{2}+1}$ then 
the integral in (\ref{wrong-3})
is $I(\tfrac{5}{3},1)$, where 
\begin{equation}
I(p,q) = \int_{0}^{\infty} \frac{dx}{\sqrt{\sigma(x,p)} \sqrt{\sigma(x,q)} \, 
\sqrt{\sigma(x,p) + \sigma(x,q)}}. 
\end{equation}
\noindent
{\em The correct value of} (\ref{wrong-3}) 
{\em has eluded us}. \\

The reader is surely aware that often typos or errors 
could have profound consequences. In a letter to Larry Glasser, 
regarding (\ref{wrong-3}), we mistyped the function $\varphi(x)$ 
of (\ref{varphi}) as
\begin{equation}
\varphi(x) = 1 + \frac{4x^{2}}{3 (1+x^{2})}. 
\end{equation}
\noindent
Larry, a consumated integrator, replied with 
$\begin{displaystyle}\sqrt{3} \left( \text{Tanh}^{-1}\sqrt{2\omega} 
- \tfrac{1}{\sqrt{2}} 
\text{Tanh}^{-1}\sqrt{\omega} \right)
\end{displaystyle}$ where
$\omega = 
\left( \sqrt{7} - \sqrt{3} \right)/2 \sqrt{7}$.   Beautiful, but it 
does not help with (\ref{wrong-3}). 
The editors of \cite{gr4} have found an alternative to this quandry: the 
latest edition \cite{gr} has no entry $3.248.5$.

Another example of errors in \cite{gr} has been discussed in the 
American Mathematical Journal
by E. Talvila \cite{talvila1}. Several entries, starting 
with $3.851.1$ \cite{gr4}
\begin{equation}
\int_{0}^{\infty} x \sin(ax^{2}) \sin(2bx) \, dx = 
\frac{b}{2a} \sqrt{\frac{\pi}{2a}} \left[ \cos \frac{b^{2}}{a} + \sin 
\frac{a^{2}}{b} \right]
\end{equation}
\noindent
are shown to be incorrect. This time, the errors are more dramatic: the 
integrals are divergent.  These entries do not appear 
in the latest edition \cite{gr}. 

The website 
\texttt{http://www.math.tulane.edu/$\sim$vhm/Table.html}
has the  goal {\em to provide proofs and context to the entries 
in} \cite{gr}. The 
example chosen for the present article is taken from Section 
$3.411$ consisting of  $32$ entries. The integrands are combinations of 
rational functions of powers and exponentials and the domain of integration
is the whole real line or the half line $(0, \infty)$. There is a single 
exception: entry $3.411.5$ states that
\begin{equation}
\int_{0}^{\ln 2} \frac{x \, dx}{1 - e^{-x}} = \frac{\pi^{2}}{12}.
\label{strange}
\end{equation}
\noindent
The next section presents an elementary proof of (\ref{strange}).

\section{A reduction argument} \label{sec-special}
\setcounter{equation}{0}

The expansion of the integrand in (\ref{strange}) as a geometric 
series yields
\begin{equation}
\frac{x}{1-e^{-x}} = x + \sum_{k=1}^{\infty} x e^{-kx}.
\end{equation}
\noindent
Term-by-term integration produces the following expressions
\begin{equation}
\int_{0}^{a} \frac{x \, dx}{1-e^{-x}} \, dx = 
\frac{1}{2}a^{2} - \sum_{k=1}^{\infty} 
\frac{e^{-ak}}{k^{2}} + \sum_{k=1}^{\infty} 
\frac{1}{k^{2}} - a \sum_{k=1}^{\infty} \frac{e^{-ak}}{k}. 
\label{whole-one}
\end{equation}
\noindent
The complexity of these three series decreases as one moves from left to
right. We now compute each term in (\ref{whole-one}), individually. 

\noindent
{\em The third series}. Integrating the geometric series 
$\begin{displaystyle} \sum_{k=0}^{\infty} x^{k} = 1/(1-x)
\end{displaystyle}$ yields $\begin{displaystyle} 
\sum_{n \geq 1} \tfrac{x^{n}}{n} = \ln(1-x) \end{displaystyle}$, which is 
valid for $|x|<1$. Evaluating at $x=e^{-a}$ gives 
\begin{equation}
\sum_{k=1}^{\infty} \frac{e^{-ak}}{k} = \ln(1- e^{-a}). 
\end{equation}

\noindent
{\em The second series}. The Riemann zeta function
\begin{equation}
\zeta(s) = \sum_{n=1}^{\infty} \frac{1}{n^{s}}
\label{zeta-def}
\end{equation}
\noindent
plays a prominent role in the evaluation of the remaining $31$ entries in 
Section $3.411$. Indeed, the first of these
\begin{equation}
\int_{0}^{\infty} \frac{x^{\nu-1} \, dx}{e^{\mu x} -1} = \frac{1}{\mu^{\nu}}
\Gamma(\nu) \zeta(\nu) 
\end{equation}
\noindent
is the classical integral representation for $\zeta(\nu)$. It 
is becoming that the special value 
\begin{equation}
\zeta(2) = \frac{\pi^{2}}{6}
\end{equation}
\noindent
appears as the second series in (\ref{whole-one}).   \\

\noindent
{\em The first series}. The second series in (\ref{whole-one}) is the 
only remaining part, we are alluding to the 
function $\begin{displaystyle} \sum_{k \geq 1} x^{k}/k^{2} \end{displaystyle}$
evaluated at $x= e^{-a}$. This is the famous {\em polylogarithm}
studied by Euler. See the introduction to \cite{lewin3} for a historical 
perspective. Aside from the series representation
\begin{equation}
\text{PolyLog}(2,x) = \sum_{k=1}^{\infty} \frac{x^{k}}{k^{2}}, 
\end{equation}
\noindent
there is a natural integral expression
\begin{equation}
\text{PolyLog}(2,x) =  - \int_{0}^{x} \frac{\ln(1-t)}{t} \, dt. 
\end{equation}

\noindent
Therefore, the identity (\ref{whole-one}) reduces to 
\begin{equation}
\int_{0}^{a} \frac{x \, dx}{1-e^{-x}} = 
\frac{1}{2}a^{2} - \text{PolyLog}[2,e^{-a}] + \frac{\pi^{2}}{6} - 
a \ln(1-e^{-a}),
\end{equation}
\noindent
and entry $3.411.5$ corresponds to the special value 
\begin{equation}
\text{PolyLog}[2, \tfrac{1}{2} ] = \frac{\pi^{2}}{12} - \frac{1}{2}\ln^{2}2.
\end{equation}
\noindent
Equivalently,
\begin{equation}
\sum_{k=1}^{\infty} \frac{1}{2^{k-1}k^{2}}  = 
\frac{\pi^{2}}{6} - \ln^{2}2.
\label{simple-polylog}
\end{equation}

Euler proved the functional equation 
\begin{equation}
\text{PolyLog}[2,x] + \text{PolyLog}[2,1-x] = \frac{\pi^{2}}{6} - 
\ln x \ln (1-x) 
\end{equation}
\noindent
for the polylogarithm function. In particular, the case $x = \tfrac{1}{2}$ 
gives (\ref{simple-polylog}).

\section{An elementary computation of the first series} \label{sec-proof}
\setcounter{equation}{0}

A series for $\ln^{2}2$ can be obtained by squaring
$\begin{displaystyle} \ln 2 = - \sum_{n \geq 1} \frac{1}{n2^{n}} 
\end{displaystyle}$ so that
\begin{eqnarray}
\ln^{2}2 & = & \left( \sum_{n=1}^{\infty} \frac{1}{n2^{n}} \right)
          \times \left( \sum_{m=1}^{\infty} \frac{1}{m2^{m}} \right) 
  =  \sum_{n,m \geq 1} \frac{1}{nm2^{n+m}} \nonumber \\
 & = & \sum_{r=1}^{\infty} \left( \sum_{m=1}^{r-1} \frac{1}{(r-m)m} 
\right) \frac{1}{2^{r}}. \nonumber
\end{eqnarray}
The partial fraction decomposition 
$\begin{displaystyle}  \frac{1}{(r-m)m} = \frac{1}{r} 
\left( \frac{1}{m} + \frac{1}{r-m} \right) 
\end{displaystyle}$ gives 
\begin{equation}
\ln^{2}2 = \sum_{r=1}^{\infty} \frac{H_{r-1}}{r2^{r-1}}, 
\end{equation}
\noindent
where $\begin{displaystyle} H_{r-1}= 1 + \tfrac{1}{2} + \cdots + \tfrac{1}{r-1}
\end{displaystyle}$ is the harmonic number. Therefore,
\begin{eqnarray}
\ln^{2}2 + \sum_{k=1}^{\infty} \frac{1}{2^{k-1}k^{2}} & = & 
\sum_{r=1}^{\infty} \frac{H_{r-1}}{r2^{r-1}} + \sum_{k=1}^{\infty} 
\frac{1}{2^{k-1}k^{2}} \nonumber \\
& = & \sum_{r=1}^{\infty} \frac{1}{r2^{r-1}} \left( H_{r-1} + \frac{1}{r} 
\right) \nonumber \\
& = & \sum_{r=1}^{\infty} \frac{H_{r}}{r2^{r-1}}. \nonumber
\end{eqnarray}
\noindent
It remains to verify that this last series is in fact $\zeta(2)$.  \\

The representation of the harmonic number as
\begin{equation}
H_{r} = \int_{0}^{1} \frac{1-x^{r}}{1-x} \, dx
\end{equation}
\noindent
gives the desired step. Indeed, if $c_{r}$ is a sequence of real numbers and 
$\alpha$ is fixed, then 
\begin{equation}
\sum_{r=1}^{\infty} c_{r}H_{r} \alpha^{r}  =  
\int_{0}^{1} \frac{1}{1-x} \sum_{r=1}^{\infty} (1-x^{r})c_{r} \alpha^{r} \, dx.
\end{equation}
\noindent
Thus the function 
\begin{equation}
f(x) = \sum_{r=1}^{\infty} c_{r}x^{r}
\end{equation}
\noindent
appears in the integral representation
\begin{equation}
\sum_{r=1}^{\infty} c_{r}H_{r} \alpha^{r} =  \int_{0}^{1} \frac{f(\alpha) - 
f(\alpha x)}{1-x} \, dx. 
\end{equation}
\noindent
In the present case, $c_{r} = 1/r2^{r-1}$ and $f(x) = -2 \ln(1-x/2)$. Therefore,
\begin{equation}
\sum_{r=1}^{\infty} \frac{H_{r}}{r2^{r-1}} = 
\int_{0}^{1} \frac{2 \ln 2 + 2 \ln(1- x/2)}{1-x} \, dx = 
2 \int_{0}^{1} \frac{\ln(1+y)}{y} \, dy. 
\end{equation}
\noindent
This last integral is computable via (\ref{simple-1}) and 
we have come full circle. 

The technique described above in exhibiting an elementary proof of 
(\ref{simple-polylog}) can be applied to
\begin{equation}
\sum_{n=1}^{\infty} \frac{H_{n-1}}{n^{2}} = \zeta(3). 
\end{equation}
\noindent
J. Borwein and D. Bradley \cite{borwein-bradley1} have given $32$ proofs 
of this charming identity. 

\medskip

\noindent
{\bf Acknowledgements.} The second author acknowledges the partial support of 
NSF-DMS 0713836.

\end{document}